\documentclass[11pt,twoside]{article}

\usepackage{amsfonts}
\usepackage[mathscr]{euscript}

\newcommand{\lbl}[1]{\label{#1}}

\newtheorem{theo}{Theorem}[section]
\newtheorem{prop}{Proposition}[section]
\newtheorem{lem}{Lemma}[section]
\newtheorem{remark}{Remark}[section]
\newtheorem{cor}{Corollary}[section]

\newcommand{\be}{\begin{equation}}
\newcommand{\ee}{\end{equation}}
\newcommand\bes{\begin{eqnarray}} \newcommand\ees{\end{eqnarray}}
\newcommand{\bess}{\begin{eqnarray*}}
\newcommand{\eess}{\end{eqnarray*}}
\newcommand\bedd{\bes\left\{\begin{array}{ll}\smallskip}
\newcommand\eedd{\end{array}\right.\ees}

\newcommand\ep{\varepsilon}
\newcommand\kk{\left}
\newcommand\rr{\right}
\newcommand\nm{\nonumber}
\newcommand\dd{\displaystyle}
\newcommand\vp{\varphi}

\newcommand\lm{\lambda}
\newcommand\ty{\infty}

\newcommand\ff{, \ \ \forall \ }

\begin{document}\thispagestyle{empty}

\begin{center}{\Large\bf The diffusive logistic equation with a free boundary}\\[2mm]
 {\Large\bf and sign-changing coefficient}\footnote{This work was supported by NSFC Grant 11371113}\\[5mm]
 {\Large  Mingxin Wang\footnote{E-mail: mxwang@hit.edu.cn; Tel: 86-15145101503; Fax: 86-451-86402528}}\\[2mm]
{Natural Science Research Center, Harbin Institute of Technology, Harbin 150080, PR China}
\end{center}

\begin{quote}
\noindent{\bf Abstract.} This short paper concerns a diffusive logistic equation with a free boundary and sign-changing coefficient, which is formulated to study the spread of an invasive species, where the free boundary represents the expanding front. A spreading-vanishing dichotomy is derived, namely the species either successfully spreads to the right-half-space as time $t\to\infty$ and survives (persists) in the new environment, or it fails to establish and will extinct in the long run. The sharp criteria for spreading and vanishing is also obtained. When spreading happens, we estimate the asymptotic spreading speed of the free boundary.


 \end{quote}

\def\theequation{\arabic{section}.\arabic{equation}}

 \section{Introduction}
 \setcounter{equation}{0} {\setlength\arraycolsep{2pt}

Understanding the nature of establishment and spread of invasive species is a central problem in invasion ecology. A lot of mathematicians have made efforts to develop various invasion models and investigated them from a viewpoint of mathematical ecology, refer to \cite{CC1}-\cite{CLL}, \cite{DG}-\cite{LL}, \cite{PZ} and \cite{Wjde}-\cite{ZX} for example. Most theoretical approaches are based on or start with single-species models. In consideration of the environmental heterogeneity, the following problem
\bess
 \left\{\begin{array}{lll}
 u_t-d\Delta u=u(m(x)-u), \ \ &t>0,\ \ x\in\Omega,\\[1mm]
 B[u]=0,  &t\ge 0, \ \ x\in\partial\Omega,\\[1mm]
 u(0,x)=u_0(x),  &x\in\Omega
 \end{array}\right.
 \eess
is a typical one to describe the spread of invasive species and has received an astonishing amount of attention, see, for example \cite{CC1, Ni} and the references therein. In this model, $u(t,x)$ represents the population density; constant $d>0$ denotes the diffusion (dispersal) rate; the function $m(x)$ accounts for the local growth rate (intrinsic growth rate) of the population and is positive on favorable habitats and negative on unfavorable ones; $\Omega$ is a bounded domain of $\mathbb R^N$; the boundary operator $B[u]=\alpha u+\beta \frac{\partial u}{\partial \nu}$, $\alpha$ and $\beta$ are non-negative functions and $\alpha+\beta>0$, $\nu$ is the outward unit normal vector of the boundary $\partial\Omega$. The corresponding systems with heterogeneous environment have also been studied extensively, please refer to \cite{CC2, CLL, LL, Ni} and the references cited therein.

To realize the spreading mechanism of an invading species (how fast spreads into new territory, and what factors influence the successful spread),  Du and Lin \cite{DLin} proposed the following free boundary problem of the diffusive logistic equation
 \bes
 \left\{\begin{array}{lll}
 u_t-du_{xx}=u(a-bu), &t>0,\ \ 0<x<h(t),\\[1mm]
 u_x(t,0)=0,\ \ u(t,h(t))=0,\ \ &t\ge 0,\\[1mm]
 h'(t)=-\mu u_x(t,h(t)),&t\ge 0,\\[1mm]
 h(0)=h_0,\ \ u(0,x)=u_0(x),&0\le x\le h_0,
 \end{array}\right.\label{1.1}
 \ees
where $x=h(t)$ is the moving boundary to be determined; $a,b,d,h_0$ and $\mu$ are
given positive constants, $h_0$ denotes the size of initial habitat, $\mu$ is the ratio of expanding speed of the free boundary and population gradient at expanding front, it can also be considered as the ``moving parameter"; $u_0$ is a given positive initial function. They have derived various interesting results.

\vskip 2pt Since then, this kind of problems describing the spread by free boundary have been studied intensively. For example, when the boundary condition $u_x=0$ at $x=0$ in (\ref{1.1}) is replaced by $u=0$, such free boundary problem was studied by Kaneko \& Yamada \cite{KY}. Du \& Guo \cite{DG, DG1}, Du, Guo \& Peng \cite{DGP} and Du \& Liang \cite{DLiang} considered the higher space dimensions, heterogeneous environment and time-periodic environment case, where the heterogeneous environment coefficients were required to have positive lower and upper bounds. Peng \& Zhao \cite{PZ} studied the seasonal succession case. Instead of $u(a-bu)$ by a general function $f(u)$, this problem has been investigated by Du \& Lou \cite{DLou} and Du, Matsuzawa \& Zhou \cite{DMZ}. The diffusive competition system with a free boundary has been studied by Guo \& Wu \cite{GW}, Du \& Lin \cite{DL2} and Wang \& Zhao \cite{WZjdde}. The diffusive prey-predator model with free boundaries has been studied by Wang \& Zhao \cite{Wjde, WZ, ZW}.

\vskip 2pt Recently, Zhou and Xiao \cite{ZX} studied the following diffusive logistic model with a free boundary in the heterogeneous environment:
 \bess
 \left\{\begin{array}{lll}
 u_t-du_{xx}=u(m(x)-u), &t>0,\  0<x<h(t),\\[1mm]
 u_x(t,0)=0,\ \ u(t,h(t))=0,\ \ &t\ge 0,\\[1mm]
 h'(t)=-\mu u_x(t,h(t)),&t\ge 0,\\[1mm]
 h(0)=h_0, \ u(0,x)=u_0(x),&0\le x\le h_0,
 \end{array}\right.
 \eess
where the initial function $u_0\in C^2([0,h_0])$, $u_0'(0)=u_0(h_0)=0$, $u_0'(h_0)<0$ and $u_0>0$ in $(0,h_0)$. In the {\it strong} heterogeneous environment, i.e,
 \vspace{-1mm}\begin{quote}
 {\bf(H1)}\, $m\in C^1([0,\infty))\cap L^\infty ([0,\infty))$ and $m$ changes sign in $(0,h_0)$,
 \vspace{-1mm}\end{quote}
Zhou and Xiao took $d$ and $\mu$ as variable parameters and derived some sufficient conditions for species spreading (resp. vanishing); While in the {\it weak} heterogeneous environment, i.e.,
  \vspace{-1mm}\begin{quote}
{\bf(H2)}\, $m\in C^1([0,\infty))$ and $0<m_1\leq m(x)\le m_2<\infty$ for all $x\geq 0$,
 \vspace{-1mm}\end{quote}
they obtained a spreading-vanishing dichotomy and a sharp criteria for spreading and vanishing. When spreading happens, they gave an estimate of the asymptotic spreading speed of the free boundary for $0<d\le d^*$ with some $d^*$.

\vskip 3pt Motivated by the above works, in this paper we consider the following problem
 \bes
 \left\{\begin{array}{lll}
 u_t-du_{xx}=u(m(x)-u), &t>0, \ 0<x<h(t),\\[1mm]
 B[u](t,0)=0, \ u(t,h(t))=0,\ \ &t \ge 0,\\[1mm]
 h'(t)=-\mu u_x(t,h(t)),&t\ge 0,\\[1mm]
 h(0)=h_0, \ u(0,x)=u_0(x),&0\le x\le h_0,
 \end{array}\right.\label{1.2}
 \ees
where, $B[u]=\alpha u-\beta u_x$, $\alpha$, $\beta\ge 0$ are constants and $\alpha+\beta=1$; the initial function $u_0(x)$ satisfies

 $\bullet$\, $u_0\in C^2([0,h_0])$ , $u_0>0$ in $(0,h_0)$, $B[u_0](0)=u_0(h_0)=0$.\\ Throughout this paper, we suppose that the function $m(x)$ satisfies

 \vspace{-2mm}\begin{quote}
{\bf(A)}\, $m\in C([0,\infty))\cap L^\infty ([0,\infty))$ and $m(x)$ is positive somewhere in $(0,\ty)$.
  \vspace{-2mm}\end{quote}
Actually, if $m(x)\le 0$ in $(0,\ty)$, the problem (\ref{1.2}) may not have the  biological background.

\vskip 3pt The objective of this paper is to study the dynamics of (\ref{1.2}) under weaker assumptions on the heterogeneous environment function $m(x)$.
In Section 2, we shall give the global existence, uniqueness, regularity and estimate of $(u,h)$. Especially, the uniform estimates of $\|u(t,\cdot)\|_{C^1[0,\,h(t)]}$ for $t\ge 1$ and $\|h'\|_{C^{\nu/2}([n+1, n+3])}$ for $n\ge 0$ are obtained directly regardless of the size of $h_\ty$, which is different from the previous works.
Section 3 is devoted to the sharp criteria for spreading and vanishing. We shall use the pairs $(h_0,\mu)$ and $(d, \mu)$, respectively, as varying parameters to describe the sharp criteria. In Section 4, we study the long time behavior of $u$ for spreading case. To this aim, in this section we first discuss the existence and uniqueness of the positive solution to a corresponding stationary problem. As a consequence of the results obtained in Sections 3 and 4, a spreading-vanishing dichotomy is obtained. In Section 5 we estimate the asymptotic spreading speed of the free boundary when spreading occurs. The last section is a brief discussion.

We remark that for the higher dimensional and radially symmetric case of (\ref{1.2}), the methods of this paper are still valid and the corresponding results can be retained. Besides, the present short paper can be regarded as the simplify, improvement and generalization of \cite{ZX} in some sense.

\section{Global existence, uniqueness and estimate of the solution $(u,h)$}
\setcounter{equation}{0}{\setlength\arraycolsep{2pt}

In this section, we give the existence, uniqueness, regularity and estimate of solution.

\begin{theo}\lbl{th2.1} \ Problem {\rm(\ref{1.2})} has a unique global solution $(u,h)$, and for some $\nu\in(0,1)$,
   \bes
  u\in C^{\frac{1+\nu} 2,1+\nu}(D_\infty),\ \ h\in C^{1+\frac{\nu}2}(0,\infty),
  \lbl{2.1}\ees
where $D_\infty=\big\{(t,x):\,t\in(0,\infty),\, x\in\big[0,h(t)\big]\big\}$. Furthermore, there exist positive constants $M=M(\|m,u_0\|_\ty)$ and  $C=C(\mu,\|m,u_0\|_\ty)$, such that
  \bes
  &0<u(t,x)\leq M, \ \ 0<h'(t)\leq \mu M, \ \ \forall \  t> 0,\ 0<x<h(t),&\lbl{2.2}\\[1mm]
  &\|h'\|_{C^{\nu/2}([n+1,n+3])}\leq C\ff n\geq 0, \ \ \
 \|u(t,\cdot)\|_{C^1([0,\,h(t)])}\leq C, \ \ \forall \  t\ge 1.&\lbl{2.3}\ees
\end{theo}

{\bf Proof.}\, Noting that the function $m$ is bounded, and applying the methods used in \cite{ChenA, DLin} with some modifications, we can prove that (\ref{1.2}) has a unique global solution $(u,h)$, and satisfies (\ref{2.1}) and the first estimate of (\ref{2.2}). The details are omitted here. Because of the condition {\bf(A)}, the regularity of $(u,h)$ can not be promoted.

\vskip 2pt Now we prove $h'(t)>0$. Firstly, as $u>0$ for $0<x<h(t)$ and $u=0$ at $x=h(t)$, we see that $u_x(t,h(t))\le 0$ and so $h'(t)\ge 0$. Since we only know $h\in C^{1+\frac\nu 2}([0,\infty))$, it can not be guaranteed that the domain $D_\infty$ has an {\it interior sphere} property at the right boundary $x=h(t)$. Hence, the Hopf boundary lemma cannot be used directly to get $h'(t)>0$. To solve this, we use a transformation to straighten the free boundary $x=h(t)$. Define
$y=x/h(t)$ and $w(t,y)=u(t,x)$. A series of detailed calculation yield
 \bess
\left\{\begin{array}{ll}
w_t-d\zeta(t)w_{yy}-\xi(t,y)w_y=w[m(h(t)y)-w], \ &t>0,\ 0<y<1,\\[1mm]
\big(\alpha w-\frac{\beta}{h(t)}w_y\big)(t,0)=0, \ \ w(t,1)=0,\ \ &t\ge 0,\\[1mm]
w(0,y)=u_0(h_0y),&0\leq y\leq 1,
\end{array}\right.
 \eess
where $\zeta(t)=h^{-2}(t)$, $\xi(t,y)=yh'(t)/h(t)$. This is an initial and boundary value problem with fixed boundary. Since $w>0$ for $t>0$ and $0<y<1$, by the Hopf boundary lemma, we have $w_y(t,1)<0$ for $t>0$. This combines with the relation $u_x=h^{-1}(t)w_y$ derives that $u_x(t,h(t))<0$, and so $h'(t)>0$ for $t>0$. The proof of $h'(t)\leq \mu M$ is similarly to that in \cite{DLin}.

Now we prove (\ref{2.3}). For the integer $n\ge 0$, let $w^n(t,y)=w(t+n,y)$, then we have
 \bess
\left\{\begin{array}{ll}
w^n_t-d\zeta(t+n)w^n_{yy}-\xi(t+n,y)w^n_y=w^n[m(h(t+n)y)-w^n], \ &t>0,\ 0<y<1,\\[1mm]
\big(\alpha w^n-\frac{\beta}{h(t+n)}w^n_y\big)(t,0)=0, \ \ w^n(t,1)=0,\ \ &t
\ge 0,\\[1mm]
w^n(0,y)=u(n,h(n)y),&0\leq y\leq 1.
\end{array}\right.
 \eess
Noticing (\ref{2.2}), apply the interior $L^p$ estimate (see \cite[Theorems 7.15 and 7.20]{Lie}) and embedding theorem, we can find a constant $C>0$ independent of $n$ such that
$\|w^n\|_{C^{\frac{1+\nu}2, 1+\nu}([1,3]\times[0,1])}\le C$ for all $n\ge 0$.
This implies $\|w\|_{C^{\frac{1+\nu}2, 1+\nu}(E_n)}\le C$, where $E_n=[n+1,n+3]\times[0,1]$.  This fact combined with
$h'(t)=-\mu u_x(t,h(t))$, $u_x(t,h(t))=h^{-1}(t)w_y(t,1)$ and $0<h'(t)\le\mu M$,
allows us to get the first estimate of (\ref{2.3}). Since these rectangles $E_n$ overlap and $C$ is independent of $n$, one has
$\|w\|_{C^{0, 1}([1,\infty)\times[0,1])}\le C$. Using $u_x=h^{-1}(t)w_y$ again, the second estimate of (\ref{2.3}) is obtained. \ \ \ \fbox{}

It follows from Theorem \ref{th2.1} that $h(t)$ is
monotonic increasing. Therefore, there exists $h_\infty\in(0,\infty]$ such that $\lim_{t\to\infty} h(t)=h_\infty$.

\section{Sharp criteria for spreading and vanishing}
\setcounter{equation}{0}{\setlength\arraycolsep{2pt}

We first prove that if $h_\infty<\infty$ then $\lim_{t\to\infty}\max_{0\leq x\leq h(t)}u(t,x)=0$. This conclusion will help us to establish the sharp criteria for spreading and vanishing.

\begin{lem}\lbl{l3.1}\, Let $d,\mu$ and $B$ be as above, $c\in\mathbb{R}$. Assume that $s\in C^1([0,\infty))$, $w\in C^{\frac{1+\nu}2,1+\nu}([0,\infty)\times[0,s(t)])$ and satisfy $s(t)>0$, $w(t,x)>0$ for $t\ge 0$ and $0<x<s(t)$. We further suppose that
$\lim_{t\to\infty} s(t)<\infty$, $\lim_{t\to\infty} s'(t)=0$ and there exists a constant $C>0$ such that $\|w(t,\cdot)\|_{C^1[0,\,s(t)]}\leq C$ for $t>1$. If $(w,s)$ satisfies
  \bess\left\{\begin{array}{lll}
 w_t-dw_{xx}\geq cw, &t>0, \ 0<x<s(t),\\[.5mm]
 B[w]=0, \ &t\ge 0, \ x=0,\\[.5mm]
 w=0, \ s'(t)\geq-\mu w_x, \ &t\ge 0, \ x=s(t),
 \end{array}\right.\eess
then $\lim_{t\to\infty}\max_{0\leq x\leq s(t)}w(t,x)=0$.
 \end{lem}

{\bf Proof}.\, When $\alpha=0$ or $\beta=0$, this is exactly \cite[Proposition 3.1]{Wjde}. When $\alpha>0$ and $\beta>0$, that proof is still valid. The details are omitted here. \ \ \ \fbox{}

\vskip 2pt Applying (\ref{2.3}) and Lemma \ref{l3.1}, we have the following theorem.

\begin{theo}\lbl{th3.1} \ Let $(u,h)$ be the solution of {\rm(\ref{1.2})}. If
$h_\infty<\infty$, then $\lim_{t\to\infty}\,\max_{0\leq x\leq h(t)}u(t,x)=0$. This shows that if the species cannot spread successfully, it will extinct in the long run.
\end{theo}

For any given $\ell>0$, let $\lm_1(\ell;d,m)$ be the first eigenvalue of
 \bes\left\{\begin{array}{ll}
 -d\phi''-m(x)\phi=\lm\phi, \ \ 0<x<\ell,\\[1mm]
 B[\phi](0)=0, \ \ \phi(\ell)=0.
 \end{array}\right.\lbl{3.2}\ees
Remember the boundary condition $\phi(\ell)=0$ and $m(x)$ is bounded, the following conclusions are well known (see, for example, \cite{CC2, Ni, Wang}).

 \begin{prop}\lbl{p3.1} \,{\rm(i)} $\lm_1(\ell;d,m)$ is continuous in $d,\,m$ and $\ell$;

{\rm(ii)} $\lm_1(\ell;d,m)$ is strictly increasing in $d$, strictly decreasing in $m$ and $\ell$;

{\rm(iii)} $\lim_{d\to\infty}\lm_1(\ell;d,m)=\lim_{\ell\to 0^+}\lm_1(\ell;d,m)=\infty$, $\lim_{d\to 0^+}\lm_1(\ell;d,m)=-\max_{[0,\ell]}m(x)$.
\end{prop}

\begin{lem}\lbl{l3.2}\, If $h_\infty<\infty$, then $\lm_1(h_\infty;d,m)\geq 0$.
\end{lem}

{\bf Proof}.\, We assume $\lm_1(h_\infty;d,m)<0$ to get a contradiction. By the continuity of $\lm_1(\ell;d,m)$ in $\ell$ and $h(t)\to h_\infty$, there exists $\tau\gg 1$ such that $\lm_1(h(\tau);d,m)<0$.
Let $w$ be the solution of
  $$\left\{\begin{array}{ll}
 w_t-dw_{xx}=w\dd\left(m(x)-w\right), \ \ &t\ge \tau, \ \, 0<x<h(\tau),\\[1mm]
 B[w](t,0)=w(t,h(\tau))=0,&t\ge \tau,\\[1mm]
  w(\tau,x)=u(\tau,x),&0\le x\le h(\tau).
  \end{array}\right.$$
Then $u\geq w$ in $[\tau,\ty)\times[0,h(\tau)]$. As $\lm_1(h(\tau);d,m)<0$, we have  $\lim_{t\to\ty}w(t,x)=z(x)$ uniformly on $[0,h(\tau)]$, where $z$ is the unique positive solution of
  \[\left\{\begin{array}{ll}
 -dz''=z\dd\left(m(x)-z\right),\ \ &0<x<h(\tau),\\[1mm]
  B[z](0)=z(h(\tau))=0.&\end{array}\right.\]
Hence, $\liminf_{t\to\infty} u(t,x)\geq z(x)>0$ in
$(0,h(T))$. This contradicts Theorem \ref{th3.1}. \ \ \ \fbox{}

\vskip 2pt The following lemma is the analogue of \cite[Lemma 3.5]{DLin} and the proof
 will be omitted.

\begin{lem} $($Comparison principle$)$\label{l3.3}\, Let $\bar h\in C^1([0,\infty))$ and $\bar h>0$ in $[0,\infty)$, $\bar u\in C^{0,1}(\overline{O})\cap C^{1,2}(O)$, with $O=\{(t,x): t>0,\, 0<x<\bar h(t)\}$. Assume that $(\bar u, \bar h)$ satisfies
 \bess\left\{\begin{array}{ll}
  \bar u_t-d\bar u_{xx}\geq\bar u(m(x)-\bar u),\ \ &t>0, \ 0<x<\bar h(t),\\[1mm]
 B[\bar u](t,0)\geq 0, \ \bar u(t,\bar h(t))=0, \ &t\ge 0,\\[1mm]
  \bar h'(t)\geq-\mu\bar u_x(t,\bar h(t)),\ \ &t\ge 0.
 \end{array}\right.\eess
If $\bar h(0)\geq h_0$, $\bar u(0,x)\geq 0$ in $[0,\bar h(0)]$, and $\bar u(0,x)\geq u_0(x)$ in $[0,h_0]$. Then the solution $(u,h)$ of {\rm(\ref{1.2})} satisfies $h(t)\leq\bar h(t)$ in $[0,\infty)$, and $u\leq\bar u$ in $D$, where $D=\{(t,x): t\geq 0,\, 0\leq x\leq h(t)\}$.
\end{lem}

\begin{lem}\label{l3.4}\, If $\lm_1(h_0;d,m)>0$, then there exists $\mu_0>0$, depending on $d,h_0,m(x)$ and $u_0(x)$, such that $h_\infty<\ty$ provided $\mu\leq\mu_0$. By Lemma {\rm\ref{l3.2}}, $\lm_1(h_\infty;d,m)\geq 0$ for $\mu\leq\mu_0$ .
\end{lem}

{\bf Proof}.\, The idea comes from \cite{DLin, GW, Wjde}, but the proof given here is more simple. Let $\phi$ be the corresponding positive eigenfunction to $\lambda_1:=\lm_1(h_0;d,m)$. Noting that $\phi'(h_0)<0$, $\phi(0)>0$ when $\beta>0$, and $\phi'(0)>0$ when $\beta=0$, it is easy to see that there exists $k>0$ such that
 \bes
 x\phi'(x)\le k\phi(x) , \ \ \forall \  \ 0\le x\le h_0.
 \lbl{3.3}\ees
Let $0<\delta,\,\sigma<1$ and $K>0$ be constants, which will be determined later. Set
 \bess
 \dd s(t)=1+2\delta-\delta {\rm e}^{-\sigma t}, \ \
 v(t,x)=K{\rm e}^{-\sigma t}\phi\left(x/s(t)\right), \ \ \ t\geq 0,\ \ 0\leq x\leq h_0s(t).
 \eess

Firstly, for any given $0<\varepsilon\ll 1$, since $m(x)$ is uniformly continuous in $[0,3h_0]$, it is easy to see that there exists $0<\delta_0(\ep)\ll 1$ such that, for all $0<\delta\le\delta_0(\ep)$ and $0<\sigma<1$,
 \bes
 \left|s^{-2}(t)m\left(x/s(t)\right)-m(x)\right|\le\varepsilon, \ \ \forall \ t>0, \ \ 0\leq x\leq h_0s(t).
\lbl{3.4}\ees

Denote $y=x/s(t)$. Owing to (\ref{3.3}), (\ref{3.4}) and $\lambda_1>0$, the direct calculation yields,
 \bes
 v_t-dv_{xx}-v(m(x)-v)&=&v\left(-\sigma+\frac {m(y)}{s^2(t)}-m(x)
-\frac{y\phi'(y)}{\phi(y)}\frac{\sigma\delta}{s(t)}{\rm e}^{-\sigma t}+\frac{\lambda_1}{s^2(t)}\right)+v^2\nm\\[1mm]
 &\geq&v(-\sigma-\varepsilon-k\sigma+\lambda_1/4)
 >0, \ \ \forall \ t>0, \ \ 0<x<h_0s(t)\qquad
 \lbl{3.5}\ees
provided $0<\sigma,\varepsilon\ll 1$. Evidently, $v(t,h_0s(t))=K{\rm e}^{-\sigma t}\phi(h_0)=0$. If either $\alpha=0$ or $\beta=0$,  then $B[v](t,0)=0$. If $\alpha,\beta>0$, then $\alpha\phi(0)=\beta\phi'(0)$ and $\phi'(0)>0$. Therefore, $B[v](t,0)=\beta K{\rm e}^{-\sigma t}\phi'(0)[1-1/s(t)]>0$ due to $s(t)>1$. In a word,
 \bes
 B[v](t,0)\ge 0, \ \ v(t,h_0s(t))=0,\ \ \forall \ t\ge 0.
 \lbl{3.6}\ees
Fix $0<\sigma,\varepsilon\ll 1$ and $0<\delta\le\delta_0(\ep)$. Thanks to the regularities of $u_0(x)$ and $\phi(x)$, we can choose a $K\gg1$ such that
 \bes
 u_0(x)\leq K\phi\left(x/(1+\delta)\right)=v(0,x),\ \ \forall \ 0\le x\le h_0.
\lbl{3.7}\ees
Thanks to $h_0s'(t)=h_0\sigma\delta{\rm e}^{-\sigma t}$ and $v_x(t,h_0s(t))=\frac 1{s(t)}K{\rm e}^{-\sigma t}\phi'(h_0)$, there exists $\mu_0>0$ such that
 \bes
 h_0s'(t)\geq-\mu v_x(t,h_0s(t)), \ \ \forall \ 0<\mu\le \mu_0, \ t\ge 0.
\lbl{3.8}\ees

Remember (\ref{3.5})-(\ref{3.8}). Applying Lemma \ref{l3.3} to $(u,h)$ and $(v,h_0s)$, it yields that $h(t)\leq h_0s(t)$ for all $t\geq 0$. Hence $h_\infty\leq h_0s(\infty)=h_0(1+2\delta)$ for all $0<\mu\leq\mu_0$. \ \ \ \fbox{}

\vskip 2pt Instead of $K$ by $\eta$, from the proof of Lemma \ref{l3.4} we see that the following lemma holds.

 \begin{lem}\lbl{l3.5} If $\lm_1(h_0;d,m)>0$, then there exist $\delta,\eta>0$,  such that $h_\infty<\ty$ provided $u_0(x)\le\eta\phi(x/(1+\delta))$ in $[0,h_0]$.
 \end{lem}

The following lemma is the analogue of \cite[Lemma 3.2]{WZjdde} and the proof will be omitted.

\begin{lem}\lbl{l3.6}\, Let $C>0$ be a constant. For any given constants $\bar h_0, H>0$, and any function $\bar u_0\in C^2([0,\bar h_0])$ satisfying $B[\bar u_0](0)=\bar u_0(\bar h_0)=0$ and $\bar u_0>0$ in $(0,\bar h_0)$, there exists $\mu^0>0$ such that when $\mu\geq\mu^0$ and $(\bar u, \bar h)$ satisfies
 \bess
 \left\{\begin{array}{ll}
   \bar u_t-d\bar u_{xx}\geq -C \bar u, \ &t>0, \ 0<x< \bar h(t),\\[1mm]
  B[\bar u](t,0)=0=\bar u(t, \bar h(t)),\ &t\geq 0,\\[1mm]
 \bar h'(t)=-\mu \bar u_x(t, \bar h(t)), \ &t\geq 0,\\[1mm]
  \bar h(0)=\bar h_0, \ \bar u(0,x)=\bar u_0(x),\ &0\leq x\leq \bar h_0,
 \end{array}\right.
 \eess
we must have $\lim_{t\to\infty}\bar h(t)>H$.
\end{lem}

To establish the sharp criteria, we define two sets.
For any given $d$, let $\sum_d=\big\{\ell>0:\,\lm_1(\ell;d,m)=0\big\}$. By the monotonicity of $\lm_1(\ell;d,m)$ in $\ell$, the set $\sum_d$ contains at most one element. For any given $\ell$, we define $\sum_\ell=\big\{d>0:\,\lm_1(\ell;d,m)=0\big\}$. Similarly, it contains at most one element.

\begin{remark}\lbl{r3.1}\,For the fixed $d>0$, due to $\lim_{\ell\to 0^+}\lm_1(\ell;d,m)=\ty$ and $\lim_{\ell\to\ty}\lm_1(\ell;d,m):=\lm_1^\ty(d,m)$ exists, we have that $\sum_d\not=\emptyset$ is equivalent to $\lm_1^\ty(d,m)<0$. As a consequence, if $m$ satisfies one of the following assumptions:
 \vspace{-1mm}\begin{quote}
 {\bf(A1)}\, There exist a constant $\rho>0$ and $y_n>x_n>0$ such that $y_n-x_n\to\ty$ as $n\to\ty$ and $m(x)\ge\rho$ in $[x_n,y_n]$;\\
  {\bf(A2)}\, There exist three constants $\rho>0$, $k>1$, $-2<\gamma\le 0$ and $x_n$ satisfying $x_n\to\ty$ as $n\to\ty$, such that $m(x)\ge\rho x^\gamma$ in $[x_n,kx_n]$.
 \vspace{-1mm}\end{quote}
Then $\lm_1^\ty(d,m)<0$, and so $\sum_d\not=\emptyset$ for all $d>0$.
\end{remark}

In fact, when the condition {\bf(A1)} holds, we use the following expression of $\lm_1(\ell;d,m)$:
 \[\lm_1(\ell;d,m)=\inf_{\phi\in H^1((0,\ell))}\frac{d\phi(0)\phi'(0)+d\int_0^\ell(\phi'(x))^2{\rm d}x-\int_0^\ell
 m(x)\phi^2(x){\rm d}x}{\int_0^\ell \phi^2(x){\rm d}x}.\]
Take a function $\phi_n$ with $\phi_n(x)=0$ in $[0, x_n]$, $\phi_n(x)=x-x_n$ in $[x_n,x_n+1]$, $\phi_n(x)=1$ in $[x_n+1,y_n-1]$ and $\phi_n(x)=y_n-x$ in $[y_n-1,y_n]$. Then  $\phi_n\in H^1((0,y_n))$, $\phi_n(0)=0$, and
 \[\int_0^{y_n}(\phi_n'(x))^2{\rm d}x=2, \ \ \int_0^{y_n}m(x)\phi_n^2(x){\rm d}x>\rho(y_n-x_n-2), \ \ \
 \int_0^{y_n}\phi_n^2(x){\rm d}x<y_n-x_n.\]
Hence, for any fixed $d>0$, we have
 \[\lm_1^\ty(d,m)<\lm_1(y_n;d,m)\le\frac{2d-\rho(y_n-x_n-2)}{y_n-x_n}\to -\rho<0 \ \ \mbox{as} \ \ n\to\ty.\]

When the condition {\bf(A2)} holds, we use the idea of {\rm\cite[Lemma 3.1]{Dong}} to derive our conclusion. Let $\lm_1(n)$ be the principal eigenvalue of
  \bess
 -d\psi''=\lm\psi, \ \ x_n<x<kx_n; \ \ \ \psi(x_n)=\psi(kx_n)=0,
 \eess
and $\psi(x)$ be the corresponding positive eigenfunction. Through a simple rescaling $\psi(x)=\Psi(x/x_n):=\Psi(y)$, we see that $\Psi(y)$ satisfies
  \bess
 -d\Psi''(y)=x_n^2\lm_1(n)\Psi(y), \ \ 1<y<k; \ \ \
 \Psi(1)=\Psi(k)=0.
 \eess
Since $\Psi>0$, we have $\lm_1^*=x_n^2\lm_1(n)$, where $\lm_1^*$ is the principal eigenvalue of
  \bess
 -d\phi''=\lm\phi, \ \ 1<x<k; \ \ \ \phi(1)=\phi(k)=0.
 \eess
Make the zero extension of $\psi$ to $[0,x_n)$, then $\psi(0)=0$ and
 \bess
 &&\dd\int_0^{kx_n}\kk[d(\psi')^2-m(x)\psi^2\rr]{\rm d}x=\dd\int_{x_n}^{kx_n}\kk[d(\psi')^2-m(x)\psi^2\rr]{\rm d}x\\[1mm]
 &=&\int_{x_n}^{kx_n}\kk[\lm_1(n)\psi^2-m(x)\psi^2\rr]{\rm d}x
\le\int_{x_n}^{kx_n}\kk(x_n^{-2}\lm_1^*-\rho k^\gamma x_n^\gamma\rr)\psi^2{\rm d}x\\[1mm]
&=&x_n^{-2}\int_{x_n}^{kx_n}\kk(\lm_1^*-\rho k^\gamma x_n^{2+\gamma}\rr)\psi^2{\rm d}x
<0 \ \ \mbox{as} \ \ n\gg 1\eess
due to $x_n\to \ty$ and $2+\gamma>0$. This implies $\lm_1(kx_n;d,m)<0$ for $n\gg 1$, and then $\lm_1^\ty(d,m)<0$.

The conditions {\bf(A1)} and {\bf(A2)} seem to be ``weaker" because $m(x)$ may be ``very negative" in the sense that both $|\{m(x)>0\}|\ll|\{m(x)<0\}|$ and $\int_0^\ty m(x){\rm d}x=-\ty$ are allowed.

\begin{remark}\lbl{r3.1a}\,For each fixed $\ell>0$, as $\lim_{d\to\ty}\lm_1(\ell;d,m)=\ty$, $\lim_{d\to 0^+}\lm_1(\ell;d,m)=-\max_{[0,\ell]}m(x)$, we see that $\sum_\ell\not=\emptyset$ is equivalent to $\max_{[0,\ell]}m(x)>0$. By the condition {\bf(A)}, we have $\max_{[0,\ell]}m(x)>0$ for each suitable large $\ell$. So, $\sum_\ell\not=\emptyset$ for such $\ell$.
\end{remark}

Now we fix $d$, and consider $h_0$ and $\mu$ as varying parameters to depict the sharp criteria for spreading and vanishing. Assume that $\sum_d\not=\emptyset$ and let $h^*=h^*(d)\in \sum_d$, i.e., $\lm_1(h^*;d,m)=0$. Recalling the estimate (\ref{2.2}), as the consequence of Lemmas \ref{l3.2}, \ref{l3.4} and \ref{l3.6}, we have

\begin{cor}\lbl{c3.1}\,{\rm(i)}\,If $h_\infty<\infty$, then $h_\infty\le h^*$. Hence, $h_0\ge h^*$ implies $h_\infty=\infty$ for all $\mu>0$;

{\rm(ii)}\,When $h_0<h^*$. There exist $\mu_0,\,\mu^0>0$, such that $h_\infty\leq h^*$ for $\mu\leq\mu_0$, $h_\infty=\infty$ for $\mu\geq\mu^0$.
\end{cor}

Finally, we give the sharp criteria for spreading and vanishing.

\begin{theo}\lbl{th3.2}\,{\rm(i)}\, If $h_0\ge h^*=h^*(d)$, then $h_\infty=\infty$ for all $\mu>0$;

{\rm(ii)}\, If $h_0<h^*$, then there exist $\mu^*>0$, depending on $d$, $m(x)$, $u_0(x)$ and $h_0$, such that $h_\infty=\infty$ for $\mu>\mu^*$, while $h_\infty\leq h^*$ for $\mu\leq\mu^*$.
\end{theo}

{\bf Proof}.\, Noticing Corollary \ref{c3.1}, by use of Lemma \ref{l3.3} and the continuity method, we can prove Theorem \ref{th3.2}. Please refer to the proof of \cite[Theorem 3.9]{DLin} for details. \ \ \ \fbox{}

\vskip 2pt When $h_0$ is fixed, $d$ and $\mu$ are regarded as the varying parameters, we have the following sharp criteria for spreading and vanishing.

\begin{theo}\lbl{th3.3}\,Assume that $\max_{[0,h_0]}m(x)>0$, and let $d^*=d^*(h_0)\in\sum_{h_0}$ {\rm(}see Remark {\rm\ref{r3.1a})}.

{\rm(i)}\, If $d\le d^*$, then $h_\infty=\infty$ for all $\mu>0$;

{\rm(ii)} If $d>d^*$ and $\sum_d\not=\emptyset$, then there exists $\mu^*>0$, depending on $d$, $m$, $u_0$ and $h_0$, such that $h_\infty=\infty$ when $\mu>\mu^*$, $h_\infty<\ty$ when $\mu\leq\mu^*$.
\end{theo}

\begin{remark}\lbl{r3.2}\, If one of {\bf(A1)} and {\bf(A2)} holds, then $\sum_d\not=\emptyset$ for any $d>0$ {\rm(}see Remark {\rm\ref{r3.1})}. \end{remark}

{\bf Proof of Theorem \ref{th3.3}}. (i)\, When $d<d^*$, we have $\lm_1(h_0;d,m)<\lm_1(h_0;d^*,m)=0$. So, $\sum_d\not=\emptyset$ and $h_0>h^*(d)$. When $d=d^*$, we have $\lm_1(h_0;d,m)=0$ and $h_0=h^*(d)$. By Theorem \ref{th3.2}(i), $h_\infty=\infty$ for all $\mu>0$.

\vskip 2pt (ii)\, For the fixed $d>d^*$, we have $\lm_1(h_0;d,m)>\lm_1(h_0;d^*,m)=0$. By Lemma \ref{l3.4}, there exists $\mu_0>0$ such that $h_\infty<\ty$ for $\mu\leq\mu_0$. On the other hand, as $\sum_d\not=\emptyset$, there exists $H\gg 1$ such that $\lm_1(H;d,m)<0$. In view of Lemma \ref{l3.6}, there exists $\mu^0>0$ such that $h_\ty>H$ provided $\mu\ge\mu^0$, which implies $\lm_1(h_\ty;d,m)<\lm_1(H;d,m)<0$. Hence, $h_\ty=\ty$ for $\mu\ge\mu^0$ by Lemma \ref{l3.2}. The remaining proof is the same as that of \cite[Theorem 3.9]{DLin}. \ \ \ \fbox{}

\vskip 4pt When $\alpha=0$ and the condition {\bf(H2)} holds, Theorem \ref{th3.3} has been given by \cite[Theorem 5.2]{ZX}.

\section{Long time behavior of $u$ for the spreading case: $h_\infty=\infty$}
\setcounter{equation}{0}

For the vanishing case: $h_\ty<\ty$, we have known $\lim_{t\to\infty}\,\max_{0\leq x\leq h(t)}u(t,x)=0$ (cf. Theorem \ref{th3.1}). In this section we study the long time behavior of $u$ for the spreading case: $h_\infty=\infty$. To this aim, we first study the existence and uniqueness of positive solution to the stationary problem:
 \bes\left\{\begin{array}{ll}
 -du''=u\big(m(x)-u\big), \ \ 0<x<\infty,\\[1mm]
 B[u](0)=0.
 \lbl{4.1}\end{array}\right.\ees

The following lemma is a special case of \cite[Proposition 2.2]{LPW}.

\begin{lem} $($Comparison principle$)$\label{l4.1}\, Let $\ell>0$, $u_1,u_2\in C^1([0,\ell))$ be positive functions in $(0,\ell)$ and satisfy
in the sense of distributions that
 \[-du_1''-m(x)u_1+u_1^2\geq 0 \geq -du_2''-m(x)u_2+u_2^2\]
and
 \[B[u_1](0)\geq 0\ge B[u_2](0), \ \ \ \limsup_{x\to \ell}(u_2^2-u_1^2)\leq 0.\]
Then $u_1\geq u_2$ in $(0,\ell)$.
\end{lem}

\begin{theo}\lbl{th4.1}\, Assume that there exist constants $-2<\gamma\le 0$ and $m_1, m_2>0$, such that
 \bes
 m_1=\liminf_{x\to\infty}\frac{m(x)}{x^\gamma}, \ \ \
 m_2=\limsup_{x\to\infty}\frac{m(x)}{x^\gamma}.
 \lbl{4.2}\ees
Then $(\ref{4.1})$ has a unique positive solution $\hat u$ and
 \bes
 m_1\le\liminf_{x\to\infty}\frac{\hat u(x)}{x^\gamma}, \ \ \
 \limsup_{x\to\infty}\frac{\hat u(x)}{x^\gamma}\le m_2.\lbl{4.3}\ees
 \end{theo}

{\bf Proof}.\, The existence of positive solution to (\ref{4.1}) can be proved as that of \cite[Lemma 7.16]{Du}. In fact, for any large $\ell>0$, in the same way as that of \cite{LPW}, we can prove that the problem
  \bess\left\{\begin{array}{ll}
 -du''=u\big(m(x)-u\big), \ \ 0<x<\ell,\\[1mm]
 B[u](0)=0, \ \ u(\ell)=\ty
 \end{array}\right.\eess
has a unique positive solution $u_\ell$ (when $\beta=0$, this conclusion is exactly \cite[Theorem 6.15]{Du}). Following the proof of \cite[Lemma 7.16]{Du} step by step (using Lemma \ref{l4.1} instead of lemma 5.6 there), we can prove that (\ref{4.1}) has at least one positive solution.

\vskip 2pt
The uniqueness of positive solution to (\ref{4.1}) and the conclusion (\ref{4.3}) can be proved by the similar way to that of \cite[Theorem 7.12]{Du} with suitable modifications. We omit the details here. Actually, proofs of the uniqueness and  (\ref{4.3}) only rely on the properties of $m$ and $u$ at infinity, have nothing to do with the condition of $u$ at $x=0$.\ \ \ \fbox{}

It is easy to see that if the condition (\ref{4.2}) holds, then the assumption {\bf(A2)} must be true. Therefore,  $\sum_d\not=\emptyset$ by Remark \ref{r3.1}.

\begin{lem}\lbl{l4.2}\, Assume that $(\ref{4.2})$ holds. Let $h^*=h^*(d)$ satisfy $\lm_1(h^*;d,m)=0$. For $\ell>h^*$, which  implies $\lm_1:=\lm_1(\ell;d,m)<0$, let $u_\ell(x)$ be the unique positive solution of
 \bes\left\{\begin{array}{ll}\smallskip
 -du''=u\big(m(x)-u\big), \ \ 0<x<\ell,\\
 B[u](0)=0, \ \ \ u(\ell)=0.
 \lbl{4.4}\end{array}\right.\ees
Then $\lim_{\ell\to\infty}u_\ell(x)=\hat u(x)$ uniformly in $[0,L]$ for any $L>0$.
 \end{lem}

{\bf Proof}.\, Let $\phi$ be the positive eigenfunction of (\ref{3.2}) corresponding to $\lambda_1$. Since $\lm_1<0$, it is easy to verify that $\ep\phi$ and $\sup_{x\ge 0}m(x)$ are the ordered lower and upper solutions to (\ref{4.4}) provided $0<\ep\ll 1$. So, the problem (\ref{4.4}) has at least one positive solution. The uniqueness of positive solution to (\ref{4.4}) is followed by Lemma \ref{l4.1}.

\vskip 2pt By Lemma \ref{l4.1}, $u_\ell\leq \hat u$ in $[0,\ell]$, and $u_\ell$ is increasing in $\ell$. Utilizing the regularity theory and compactness argument, it follows that there exists a positive function $u$, such that $u_\ell\to u$
in $C^2_{\rm loc}([0,\infty))$ as $\ell\to\infty$, and $u$ solves (\ref{4.1}). By the uniqueness, $u=\hat u$. \ \ \ \fbox{}

\vskip 2pt Finally, we give the main result of this section.

\begin{theo}\lbl{th4.2}\,Let $(\ref{4.2})$ hold. If $h_\infty=\infty$, then $\lim_{t\to\infty}u(t,x)=\hat u(x)$ in $C_{\rm loc}([0,\ty))$.
\end{theo}

{\bf Proof}.\, Choose $K>1$ such that $K\hat u\ge u_0$ in $[0,h_0]$. Then $\vp:=K\hat u$ satisfies $\vp_t-d\vp_{xx}>\vp(m(x)-\vp)$. Let $w$ be the solution of
 \bess\left\{\begin{array}{lll}
 w_t-dw_{xx}=w(m(x)-w), &t>0,\ \ 0<x<\infty,\\[1mm]
 B[w](t,0)=0,\ \ \ &t>0,\\[1mm]
  w(0,x)=K\hat u(x), &x\geq 0.
  \end{array}\right.\eess
Then $u\leq w$, and $w$ is monotone decreasing in $t$.
Because $\hat u$ is the unique positive solution of (\ref{4.1}), by the standard method we can prove that $\lim_{t\to\infty}w(t,x)=\hat u(x)$ uniformly in $[0,L]$ for any $L>0$. As $h_\infty=\infty$, it follows that $\limsup_{t\to\infty}u(t,x)\le\hat u(x)$ uniformly in $[0,L]$.

\vskip 3pt Let $h^*=h^*(d)$ be such that $\lm_1(h^*;d,m)=0$. When $\ell>h^*$, we have $\lm_1:=\lm_1(\ell;d,m)<0$. As $h_\ty=\ty$, there exists $T\gg 1$ such that $h(t)>\ell$ for all $t\ge T$. Let $\phi$ be the positive eigenfunction of (\ref{3.2}) corresponding to $\lambda_1$. Choose $0<\sigma\ll 1$ such that $u(T, x)\geq\sigma\phi(x)$ in $[0, \ell]$ and $\sigma\phi$ is a lower solution of (\ref{4.4}). Let $u^\ell$ be the unique solution of
 \bess\left\{\begin{array}{lll}
 u_t-du_{xx}=u(m(x)-u),\ \ &t\ge T, \ \ 0<x<\ell,\\[1mm]
 B[u](t,0)=0, \ \ u(t,\ell)=0,\ \ \ &t\ge T,\\[1mm]
 u(T,x)=\sigma\phi(x),&x\in [0,\ell].
 \end{array}\right.\eess
Then $u\geq u^\ell$ in $[T,\ty)\times[0,\ell]$, and $u^\ell$ is increasing in $t$. So,  $\lim_{t\to\infty}u^\ell(t,x)=u_\ell(x)$ uniformly in $[0,\ell]$ since $u_\ell$ is the unique positive solution of (\ref{4.4}). Hence, $\liminf_{t\to\infty}u(t,x)\geq u_\ell(x)$ uniformly in $[0,\ell]$. By Lemma \ref{l4.2}, $\liminf_{t\to\infty}u(t,x)\geq\hat u(x)$ uniformly in $[0,L]$ for any $L>0$. \ \ \ \fbox{}

\vskip 3pt Here we remark that, when $\alpha=0$, Theorem \ref{th4.2} has been obtained by {\rm\cite{ZX}} under one of the following assumptions:

{\rm(i)}\, the condition {\bf(H2)} holds {\rm(}see {\rm\cite[Lemma 5.2]{ZX})};

{\rm(ii)}\, the function $m\in C^1([0,\infty))$, is positive somewhere in $(0,h_0)$ and satisfies $(\ref{4.2})$ with $\gamma=0$. The diffusion rate $d$ satisfies
$0<d\le d^*$ for some $d^*>0$ {\rm(}see {\rm\cite[Lemma 6.2]{ZX})}.

Obviously, {\bf(H2)} implies $(\ref{4.2})$ with $\gamma=0$.

\vskip 3pt
Combining Theorems \ref{th3.1}, \ref{th3.2}, \ref{th3.3} and \ref{th4.2}, we have  the following two theorems concerning spreading-vanishing dichotomy and sharp criteria for spreading and vanishing.

\begin{theo}\lbl{th4.3} Let {\rm(\ref{4.2})} hold, $d>0$ be fixed and $h^*=h^*(d)$ satisfy $\lm_1(h^*;d,m)=0$. Then either

{\rm(i)} Spreading: $h_\ty=\ty$ and $\lim_{t\to\ty}u(t,x)=\hat u(x)$ uniformly in $[0,L]$ for any $L>0$; or

{\rm(ii)} Vanishing: $h_\ty\le h^*$ and $\lim_{t\to\infty}\,\max_{0\leq x\leq h(t)}u(t,x)=0$, where $\hat u(x)$ is the unique positive solution of {\rm(\ref{4.1})}.

Moreover,

{\rm(iii)} If $h_0\ge h^*$, then $h_\infty=\infty$ for all $\mu>0$;

{\rm(iv)}\, If $h_0<h^*$, then there exist $\mu^*>0$, depending on $d$, $m$, $u_0$ and $h_0$, such that $h_\infty=\infty$ for $\mu>\mu^*$, while $h_\infty\leq h^*$ for $\mu\leq\mu^*$.
\end{theo}

\begin{theo}\lbl{th4.4} Assume that {\rm(\ref{4.2})} holds, $h_0>0$ is fixed and  $\max_{[0,h_0]}m(x)>0$. Let $d^*=d^*(h_0)\in\sum_{h_0}$  Then either

{\rm(i)} Spreading: $h_\ty=\ty$ and $\lim_{t\to\ty}u(t,x)=\hat u(x)$ uniformly in $[0,L]$ for any $L>0$; or

{\rm(ii)} Vanishing: $h_\ty<\ty$ and $\lim_{t\to\infty}\,\max_{0\leq x\leq h(t)}u(t,x)=0$.

Moreover,

{\rm(iii)} If $d\le d^*$, then $h_\infty=\infty$ for all $\mu>0$;

{\rm(iv)}\, If $d>d^*$, then there exist $\mu^*>0$, depending on $d$, $m$, $u_0$ and $h_0$, such that $h_\infty=\infty$ for $\mu>\mu^*$, while $h_\infty<\ty$ for $\mu\leq\mu^*$.
\end{theo}

\section{Asymptotic spreading speed}
\setcounter{equation}{0}{\setlength\arraycolsep{2pt}

In this section, we shall estimate the asymptotic spreading speed of the free boundary $h(t)$ when spreading occurs. Throughout this section, we assume that $(\ref{4.2})$ holds with $\gamma=0$, which implies $\sum_d\not=\emptyset$ for all $d>0$.

\vskip 2pt Let us first state a known result, which plays an important role in later discussion.

\begin{prop}\lbl{p5.1}{\rm(\cite[Proposition 4.1]{DLin})}\, Let $d$ and $c$ be given positive constants. Then for any $k\ge 0$, the problem
\bess\left\{\begin{array}{ll}
 -dw''+kw'=w(c-w), \ \ 0<x<\infty,\\[1mm]
 w(0)=0, \ \ w(\infty)=c
  \end{array}\right.\eess
has a unique positive solution $w_k(x)$.  Moreover, for each $\mu>0$, there exists a unique $k_0=k_0(\mu,c)>0$ such that $\mu w_{k_0}'(0)=k_0$.
 \end{prop}

\begin{theo}\lbl{th5.1} When $h_\infty=\infty$, we have {\rm(}no other restrictions on $d,h_0,m$ and $u_0${\rm)}
 \bes
  k_0(\mu,m_1)\le\liminf_{t\to\infty}\frac{h(t)}t, \ \ \ \limsup_{t\to\infty}\frac{h(t)}t\le k_0(\mu,m_2).
\lbl{5.1}\ees
 \end{theo}

{\bf Proof}.\,The proof is similar to those of \cite[Theorem 4.2]{DLin}, \cite[Theorem 3.6]{DG} and \cite[Theorem 6.1]{ZX}. Here we give the sketch for completeness and readers' convenience.
 
For any given $0<\varepsilon\ll 1$, by (\ref{4.2}) and (\ref{4.3}) with $\gamma=0$, there exists $\ell=\ell(\varepsilon)\gg 1$ such that
 \[m_1-\varepsilon<m(x)<m_2+\varepsilon, \ \ m_1-\varepsilon<\hat u(x)<m_2+\varepsilon, \ \ \forall \  x\ge \ell.\]
Take advantage of $h_\infty=\infty$ and Theorem \ref{th4.2}, there exists $T=T(\ell)\gg 1$ such that
 \[h(T)>2\ell, \ \ \ m_1-2\varepsilon<u(t+T, \ell)<m_2+2\varepsilon, \ \ \forall \  t>0.\]
Follow the proof of \cite[Theorem 3.6]{DG} or \cite[Theorem 6.1]{ZX} step by step, we can get (\ref{5.1}). The details are omitted here. \ \ \ \fbox{}

When $\alpha=0$, Theorem \ref{th5.1} has been given in \cite{ZX} for the case that $0<d\le d^*$ with  some $d^*>0$.
\section{Conclusion}

From the above discussions we have seen that $\lm_1^\ty(d,m):=\lim_{\ell\to\ty}\lm_1(\ell;d,m)<0$ is an essential condition. This number is only characterized by $d$ and $m$, and is independent of the moving parameter $\mu$ and initial value $u_0(x)$. It seems that $\lm_1^\ty(d,m)$ is determined by $d$ and $\int_0^\ty m(x){\rm d}x$.

The main conclusions of this paper can be briefly summarized as follows:

(I)\, If one of the following holds:

(i)\, $d$ is suitable small ($h_0$ and $m(x)$ are fixed, $m(x)$ is positive somewhere in $(0,h_0)$),

(ii)\, $m(x)$ is suitable ``larger" in the sense of ``distribution" ($h_0$ and $d$ are fixed),

(iii)\, $h_0$ is suitable ``larger" ($d$ and $m(x)$ are fixed, $m(x)$ satisfies either {\bf(A1)} or {\bf(A2)}),\\
then the species will successfully spread and survive in the new environment (maintain a positive density distribution), regardless of initial population size and value of the moving parameter.

\vskip 2pt (II)\, When the above situations are not appeared, we can control the moving parameter $\mu$ and find a critical value $\mu^*$ such that the species will spread successfully when $\mu>\mu^*$, the species fails to establish and will extinct in the long run when $\mu\le\mu^*$. The better way to reduce the moving parameter might be by controlling the surrounding environment.

These theoretical results may be helpful in the prediction and prevention of biological invasions.

\bibliographystyle{elsarticle-num}

\end{document}